\documentclass[	a4paper, 
								11pt]{article}
								
\usepackage{amssymb,amsmath,amsthm}
\usepackage{graphicx}
\usepackage{subfigure}
\usepackage{overpic}
\usepackage[table]{xcolor}
\usepackage{chngcntr}
\usepackage{lscape}
\usepackage{rotating}
\usepackage{booktabs}
\usepackage{units}

\usepackage{epstopdf}

\usepackage{multirow}
\usepackage{calc}

\usepackage{tikz}
\usepackage{tikz-cd}
\usepackage{pgfplots}
\usepackage{pgfplotstable}
\usetikzlibrary{3d,calc}
\usetikzlibrary{shapes, spy, scopes, matrix, datavisualization.formats.functions, arrows, datavisualization, plothandlers, plotmarks, svg.path, positioning}
\usetikzlibrary{external}
\tikzset{external/only named=true}
\pgfplotsset{compat=newest}
\pgfplotsset{plot coordinates/math parser=false}
\pgfplotsset{scaled y ticks = false, tick label style={/pgf/number format/fixed}}
\pgfplotsset{scaled x ticks = false, tick label style={/pgf/number format/fixed}}

\usepackage{float}

\usepackage[bookmarksopen, pdfstartview={FitV}, linkbordercolor=white]{hyperref}

\usepackage{nomencl}

\usepackage[]{algorithm2e}

\usepackage{todonotes}

\usepackage[capitalize,nameinlink]{cleveref}
	
\newtheoremstyle{mystyle}%
  {}%
  {}%
  {\itshape}%
  {}%
  {\bfseries}%
  {.}%
  { }%
  {}%
	
\theoremstyle{mystyle}
\newtheorem{theorem}{Theorem}[section]
\theoremstyle{plain}
\newtheorem{definition}[theorem]{Definition}
\theoremstyle{plain}
\newtheorem{lemma}[theorem]{Lemma}
\theoremstyle{plain}
\newtheorem{remark}[theorem]{Remark}
\theoremstyle{plain}

\theoremstyle{plain}
\newtheorem{corollary}[theorem]{Corollary}
\theoremstyle{plain}

\theoremstyle{plain}

\newcommand{\BIGOP}[1]{\mathop{\mathchoice%
{\raise-0.22em\hbox{\huge $#1$}}%
{\raise-0.05em\hbox{\Large $#1$}}{\hbox{\large $#1$}}{#1}}}

\newcommand{\BIGboxplus}{\mathop{\mathchoice%
{\raise-0.35em\hbox{\huge $\boxplus$}}%
{\raise-0.15em\hbox{\Large $\boxplus$}}{\hbox{\large $\boxplus$}}{\boxplus}}}

\makeatletter
\def\tagform@#1{\maketag@@@{\ignorespaces#1\unskip\@@italiccorr}}
\let\orgtheequation\theequation
\def\theequation{(\orgtheequation)}
\makeatother

\definecolor{UBTgreen}{RGB}{0, 146, 96}

\definecolor{UBTorangeComplementary}{RGB}{218,71,0}

\definecolor{UBTblueAdjacent}{RGB}{10,77,141}
\definecolor{UBTgreenAdjacent}{RGB}{75,191,0}

\definecolor{UBTorangeTriad}{RGB}{218,129,0}
\definecolor{UBTredTriad}{RGB}{199,0,50}

\definecolor{UBTblueTetrad}{RGB}{10,77,141}
\definecolor{UBTorangeTetrad}{RGB}{218,129,0}
\definecolor{UBTredTetrad}{RGB}{218,71,0}

\def\N{\mathbb N}	\def\Z{\mathbb Z}	\def\Q{\mathbb Q}	 \def\R{\mathbb R}

\def\H{\mathbb{H}}
\def\R{\mathbb{R}}
\def\N{\mathbb{N}}
\def\M{\mathbb{M}}

\def\V{\mathbb{V}}
\def\Y{\mathbb{Y}}

\def\cA{\mathcal{A}}

\def\cH{\mathcal{H}}
\def\cI{\mathcal{I}}
\def\cJ{\mathcal{J}}

\def\cS{\mathcal{S}}

\def\ee{\boldsymbol{e}}
\def\ff{\boldsymbol{f}}

\def\kk{\boldsymbol{k}}

\def\nn{\boldsymbol{n}}

\def\rr{\boldsymbol{r}}

\def\uu{\boldsymbol{u}}

\def\xx{\boldsymbol{x}}
\def\yy{\boldsymbol{y}}

\def\AA{\boldsymbol{A}}

\def\NN{\boldsymbol{N}}

\def\SS{\boldsymbol{S}}

\def\11{\mathbf{1}}

\def\00{\boldsymbol{0}}

\def\aalpha{\boldsymbol{\alpha}}
\def\bbeta{\boldsymbol{\beta}}

\def\llambda{\boldsymbol{\lambda}}

\def\xxi{\boldsymbol{\xi}}

\def\oomega{\boldsymbol{\omega}}

\def\XXi{\boldsymbol{\Xi}}

\def\fraku{\mathfrak{u}}

\def\transpose{\text{T}}

\usepackage{comment}

\theoremstyle{plain}

\usepackage[backend=biber,
hyperref,
backref,
bibencoding=utf8,
defernumbers=true,
style=numeric-comp,
giveninits=true,
maxnames=99,
doi=false,
isbn=false, url=false]{biblatex}
\renewbibmacro{in:}{\ifentrytype{article}{}
			{\printtext{\bibstring{in}\intitlepunct}}}
\usepackage{csquotes}
\makeatletter

\let\blx@rerun@biber\relax

\makeatother

\usepackage[toc]{appendix}

\usepackage{authblk}

\title{Numerical Aspects of the Tensor Product Multilevel Method for 
High-dimensional, Kernel-based Reconstruction on Sparse Grids}
\date{\today}

\author[1]{Markus B\"uttner}
\author[2]{R\"udiger Kempf}
\author[3]{Holger Wendland}
\affil[1]{Chair of Scientific Computing\\ 
	   Department of Mathematics\\ 
	   University of Bayreuth\\
	   95440 Bayreuth\\
	   Germany}
\affil[2,3]{Applied and Numerical Analysis\\ 
	   Department of Mathematics\\ 
	   University of Bayreuth\\ 
	   95440 Bayreuth\\
	   Germany}

\begin{document}

\maketitle

\abstract{This paper investigates the approximation of functions with finite smoothness defined on domains with a Cartesian product structure. The recently proposed \emph{tensor product multilevel method (TPML)} combines Smolyak's sparse grid method with a kernel-based residual correction technique. The contributions of this paper are twofold. First, we present two improvements on the TPML that reduce the computational cost of point evaluations compared to a naive implementation. Second, we provide numerical examples that demonstrate the effectiveness and innovation of the TPML.}



\section{Introduction}

Many current problems in physics, engineering, and even business science are higher-dimensional in nature. This could be because the problem is modeled as a partial differential equation in space-time, which results in a two/three- plus-one-dimensional problem, or because it involves uncertainty, which is modeled by a stochastic process and then parameterized, done so in uncertainty quantification. The number of parameters then corresponds to the dimension of the problem. Classical numerical methods are often numerically not feasible, even in only four-dimensional examples, because of the curse of dimensionality. One way to alleviate this curse is by means of Smolyak's method \cite{smolyak:Smolyak}. This method is tailored for problems whose domain can be decomposed into the Cartesian product of lower-dimensional domains, for example space and time, and allows us to combine numerical methods used on these lower-dimensional domains to an approximation process in high dimensions and offers a very flexible ansatz with applications in quadrature \cite{smolyak:Smolyak,gerstner:NumericalIntegration,gerstner:DomensionAdaptiveTensorProductQuadrature,hajiAli:NovelResultsAnisotropic}, interpolation \cite{conrad:AdaptiveSmolyakPseudospectralApproximations,barthelmann:HighDimensionalPolynomialInterpolationOnSparseGrids,kempf:TPML}, solving PDEs \cite{bungartz:ANoteOnTheComplexityOfSolvingPoissonEquation,nobile:AnAnisotropocSparseGrid,griebel:SparseGridsSchroedingerEquation}, machine learning \cite{garcke:DimensionAdaptiveSparseGridMachineLearning} and even financial mathematics \cite{li:SparseGridInterpolationAmericanOptions}. The list given here is by no means exhaustive. For a more complete list, we refer to \cite{garcke:SparseGridsInANutshell,griebel:SparseGrids}. 

 Recently, a new method, called the \emph{tensor product multilevel method (TPML)}, was proposed and studied in \cite{kempf:TPML}, combining the method of Smolyak with a kernel-based multilevel scheme \cite{wendland:multilevel}. This method allows us to approximate moderately high-dimensional functions defined on the Cartesian product of lower-dimensional domains using scattered data. The ability to use general lower-dimensional domains instead of only intervals is one of the exceptional features of this method. Other known methods are either spline or polynomial based and are therefore based on intervals and thus fail at solving these kinds of problems. There have been previous attempts to combine Smolyak's method with a kernel-based approach, see, e.g., \cite{levesley:MultilevelSparseKernelInterpolation}, but concentrated on non-compactly supported basis functions and did not provide a rigorous error analysis. In contrast, the paper \cite{kempf:TPML} provides a thorough error analysis for compactly supported kernels of finite smoothness. There, an explicit representation of the approximation operator is also given. However, it turns out that a naive implementation of the derived representation is not feasible because the computational cost is too large. Hence, in the present paper, we propose two different reformulations of the TPML, with the result that the method can actually be applied to real-world examples. The first reformulation strives to relocate expensive computations into an offline phase such that the numerical cost of several point evaluations of the approximation becomes cheaper compared to the naive implementation. The second reformulation is tailored to problems with nested sets of sites, which occur naturally in typically applications. There, we use the nestedness to derive a representation of the TPML that uses every data value only once, reducing the numerical costs this way.

The outline of this paper is as follows. In \cref{sec:TMPL} we briefly recall the building blocks of the TPML, i.e., Smolyak's method and the formulation of the kernel-based multilevel method using a (modified) Lagrange basis, before giving the definition of the TPML operator. In \cref{sec:ModificationsOfTPML} we introduce the two variations of the formulation of TPML. And finally in \cref{sec:NumericalExamples} we demonstrate the power of the new formulations by applying the TPML to a real-world application: the tidal flow at the Bight of Abaco, an example of an interpolation of simulation data of a space-time differential equation, and to a benchmarking problem: a cantilever beam with seven parameters. To keep \cref{sec:ModificationsOfTPML} brief, we give the proofs of the results in \cref{app:Proofs,app:Proofs2}.

\section{The Tensor Product Multilevel Method}\label{sec:TMPL}

We briefly repeat the basic ideas of the kernel-based tensor product multilevel method. For details, we refer to \cite{kempf:TPML}. The method essentially combines the kernel-based multilevel method for scattered data approximation on arbitrary low-dimensional domains with Smolyak's method, which allows us to construct an approximation method for high-dimensions with essentially the same properties as the building-block, low-dimensional methods.

\subsection{Smolyak's Method}

Smolyak's construction and subsequent generalizations, as well as the numerical approximation power and computational effort depend on a pre-determined 
\textit{index set}. This set describes how many levels we will use in each 
direction. In order to have the maximum freedom and to adapt to possible
different smoothnesses in different directions it became customary to use \textit{anisotropic index sets}, see, e.g., \cite{nobile:AnAnisotropocSparseGrid,
conrad:AdaptiveSmolyakPseudospectralApproximations,
hajiAli:NovelResultsAnisotropic,
wendlandRieger:AnisotropicSamplingInequalities}.

\begin{definition}
 For $ d \in \N $ let $ \oomega \in \R^d_+ $ and $ \ell \in \N $. The 
 \emph{anisotropic index set} $ \cI_{\oomega}(\ell,d) \subset \N^d $ is 
 defined by 
 \begin{align}\label{eq:defAnisotropicIndexSet}
  \cI_{\oomega}(\ell,d) := \left\{ \llambda \in \N^d \; : \; \sum_{j=1}
  ^d (\lambda_j - 1) \omega_j \leq \ell \omega_{min} \right\},
 \end{align}
  where $ \omega_{min} := \min_{1 \leq j \leq d} \omega_j $.

 Its \emph{surface} $ \cJ_{\oomega}(\ell,d) $ is defined by
 \begin{align}\label{eq:defAnisotropicIndexSetSurface}
  \cJ_{\oomega}(\ell,d) = \cI_{\oomega}(\ell,d) \setminus \cI_{\oomega}
  \left(\ell - \frac{\| \oomega \|_1}{\omega_{min}}, d \right).
 \end{align}
\end{definition}

The positive weight vector $ \oomega $ reflects how important each
direction is. The larger the quotient $ \omega_j / \omega_{min} $ is chosen the 
less important direction $ j $ is for us and $ \cI_{\oomega}(\ell, d) $ is less 
extended in this direction. For fixed $ \oomega $ the threshold $ \ell $ determines the maximum number of levels in each direction. The largest level in 
each direction can be directly computed by 
\begin{align*}
 \lambda_{j, max} = \left \lfloor \frac{\ell \omega_{min}}{\omega_j} \right
 \rfloor + 1.
\end{align*}

\begin{remark}
For the error analysis in \cite{kempf:TPML}, we have to assume that the weight
vector $ \oomega $ is ascendingly ordered, i.e., $ \omega_1 \leq \omega_2 \leq 
\cdots \leq \omega_d $.
\end{remark}

Given sequences of operators $ A^{(j)}_i $, $ 1 \leq j \leq d $, $ 1 \leq i \leq \lambda_{j,max} $, and setting $ A^{(j)}_0 = 0 $, the Smolyak operator in combination technique representation is given as 
\begin{align*}
    \cS_{\cI_{\oomega}(\ell,d)} := \sum_{\llambda \in \cJ_{\oomega}(\ell,d)} \sum_{\substack{\bbeta \in \{0,1\}^d \\ \llambda + \bbeta \in \cI_{\oomega}(\ell,d)}} (-1)^{|\bbeta|} \left( A^{(1)}_{\lambda_1} \otimes \cdots \otimes A^{(d)}_{\lambda_d} \right).
\end{align*}

Smolyak's method is deeply connected to \emph{sparse grids}, which we define now. To families of sites $ X_{\lambda_j}^{(j)} $, $ 1 \leq \lambda_j \leq \lambda_{j,max} $, the sparse grid $ \cH_{\cI_{\oomega}(\ell,d)} \subseteq \Omega $ is defined as
\begin{align*}
    \cH_{\cI_{\oomega}(\ell,d)} := \bigcup_{\llambda \in \cI_{\oomega}(\ell,d)} X^{(1)}_{\lambda_1} \otimes \cdots \otimes X^{(d)}_{\lambda_d}.
\end{align*}

\subsection{Kernel Multilevel Method using Lagrange Functions}

The standard kernel-based multilevel method, described in, e.g., \cite{wendland:multilevel}, is a residual correction scheme for scattered data approximation. The key ingredients are a sequence of, not necessarily nested, sets of sites \( X_1, \dots, X_L \), level-dependently scaled kernels $ \Phi_i $ with compact support, and associated \emph{local approximation spaces}
\begin{align*}
 W_i = \operatorname{span}\left\{ \Phi_i ( \cdot - \xx_{i,k})
 \; : \; \xx_{i,k} \in X_i \right \}, \quad 1 \leq i \leq L.
\end{align*}
In the standard formulation of the kernel-based multilevel method, there is only an implicit dependence of the data. Hence, the authors of \cite{kempf:TPML} found a new representation of the method in terms of \emph{Lagrange} or \emph{cardinal functions} $ \{ \chi_{i,k} \} $, that made the dependence on the data explicit. 
The defining property of Lagrange functions is that they satisfy $ \chi_{i,k} (\xx_{i,m}) = \delta_{km} $ for $ \xx_{i,m} \in X_i $.

The Lagrange basis can be computed by solving a linear system. However, the 
matrix differs whether we use interpolation or penalized least squares 
approximation. We remark that in the case of penalized least squares, the functions do not satisfy the Lagrange condition, but lead to good numerical results anyway.

\begin{definition}
For $ L \in \N $, $ 1 \leq i \leq L $ and $ 1 \leq k \leq 
N_i $ let 
\begin{align}\label{eq:definitionRVector}
 \rr_i := \left( \Phi_i( \cdot - \xx_{i,1}), \dots, 
\Phi_i( \cdot - \xx_{i,N_i}) \right)^{\transpose} \in \R^{N_i}.
 \end{align}
The $k$-th Lagrange function on level $ i $, $ \chi_{i,k} $, can then be expressed as
\begin{align}\label{eq:LagrangeFunctionVectorVectorProduct}
 \chi_{i,k} = \left(\aalpha_{i,k}\right)^{\transpose} \rr_i.
\end{align}
The coefficient vector $ \aalpha_{i} \in \R^{N_i} $ is given as
 \begin{align}\label{eq:coefficientVector}
  \aalpha_{i,k} := \begin{cases} 
				  M_i^{-1} \ee_k, &\text{for interpolation,} \\
                  (M_i + \lambda_i I)^{-1} \ee_k, &\text{for penalized 
                  least-squares approximation}
                 \end{cases}
 \end{align}
where $\left(M_i\right)_{k,m} := \left(\Phi_i( \xx_{i,k} - 
 \xx_{i,m}) \right)$ is the Gramian, $ \ee_k $ denotes the $ k $-th unit vector and 
 $ I \in \R^{N_i \times N_i} $ is the unit matrix and $ \lambda_i $ are chooseable regularization parameter.
\end{definition}

We use the following notation throughout the rest of this paper.

\begin{definition}
\begin{enumerate}
 \item For $ L \in \N $ we call a subset $ \fraku := \{ u_1, \dots, u_m \} 
 \subseteq \{1, \dots, L \} $ with $ m = \# \fraku $ elements \emph{ordered}, if $ u_1 < u_2 < \cdots < u_m $.
 \item  For an ordered set $ \fraku = \{ u_1, \dots, u_m \} $, we take 
 $ \NN_{\fraku} $ as $ \NN_{\fraku} := (N_{u_1}, \dots, N_{u_m})^{\transpose} \in \N^{\# \fraku} $.
 \item For an ordered set $ \fraku = \{ u_1, \dots, u_m \} $ we abbreviate
 \begin{align*}
  \sum_{\kk \leq \NN_{\fraku}} := \sum_{k_1 \leq N_{u_1}} \cdots \sum_{k_m \leq
  N_{u_m}},
 \end{align*}
 where we implicitly assume that the dimension of $ \kk $ fits the dimension of 
 $ \NN_{\fraku} $.
\end{enumerate}
\end{definition}

Then, the multilevel approximation operator can be expressed in the following way, see \cite[Theorem 3.9]{kempf:TPML}.

\begin{theorem}\label{thrm:RepresentationMultilevelOperator}
 Let $ \fraku = \{u_1, \dots, u_{\# \fraku} \} $ be an ordered set. Define the \emph{combined operator} $ \cI_{\fraku}: C(\Omega) \to W_{\# \fraku} $ by
 \begin{align}\label{eq:defCombinedOperator}
     \cI_{\fraku} f := \sum_{\kk \leq \NN_{\fraku}} a(\fraku, \kk) f(\xx_{u_1,k_1}) \chi_{u_{\# \fraku}, k_{\# \fraku}},
 \end{align}
 where the coefficients are given by $ a(\fraku, \kk) = 1 $, if $ \# \fraku = 1 $ and 
 \begin{align*}
     a(\fraku,\kk) = \prod_{\ell = 1}^{\# \fraku - 1} \chi_{u_{\ell}, k_{\ell}} ( \xx_{u_{\ell +1},k_{\ell+1}})
 \end{align*}
 for $ \# \fraku > 1 $.
 
 Then the multilevel approximation operator $ A_L: C(\Omega) \to \bigoplus_{i=1}^L W_i $ at level $ L$  has the representation
 \begin{align*}
     A_L(f) = \sum_{\substack{\fraku \subseteq \{1, \dots, L \} \\ 1 \leq \# \fraku \leq L}} (-1)^{\# \fraku + 1} \sum_{\kk \leq \NN_{\fraku}} a(\fraku, \kk) f(\xx_{u_1,k_1}) \chi_{u_{\# \fraku}, k_{\# \fraku}}
 \end{align*}
 for $ f \in C(\Omega) $.
\end{theorem}

\subsection{The Tensor Product Multilevel Method}

The combination of Smolyak's method and the kernel-based multilevel method yields the tensor product multilevel approximation operator for continuous functions on $ \Omega^{\otimes} = \Omega^{(1)} \times \cdots \times \Omega^{(d)}$, where $ \Omega^{(j)} \subseteq \R^{n_j} $, $ 1 \leq j \leq d $ and $ n_j \in \N $. In particular, we allow $ n_j > 1 $.

\begin{definition}\label{def:TPMLGeneralRepresentation}
 We define the \emph{tensor product multilevel approximation} to $ \ff \in C(\Omega^{\otimes}) $, 
 evaluated at $ \xx = ( \xx^{(1)}, \dots, \xx^{(d)}) \in \Omega^{\otimes} $, as 
   \begin{align}\label{eq:TPMLGeneralRepresentation}
 \begin{aligned}
  &\cA_{\cI_{\oomega}(\ell,d)}(\ff)(\xx) = 
  \sum_{\llambda \in \cJ_{\oomega}(\ell,d)} \sum_
  {\substack{\bbeta \in \{0,1\}^d \\ \llambda + \bbeta \in \cI_{\oomega}
  (\ell,d)}} \\
  &\sum_{\substack{\fraku^{(1)} \subseteq \{1, \dots, \lambda_1\} \\ 1 \leq \# 
  \fraku^{(1)} \leq \lambda_1}} \cdots \sum_{\substack{\fraku^{(d)} \subseteq 
  \{1, \dots, \lambda_d\} \\ 1 \leq \# \fraku^{(d)} \leq \lambda_d}} 
  c_{\bbeta}( \fraku^{(1)}, \dots, \fraku^{(d)}) \\ 
  &\sum_{\kk^{(1)} \leq 
  \NN_{\fraku^{(1)}}} \cdots \sum_{\kk^{(d)} \leq \NN_{\fraku^{(d)}}} 
  \ff(\xx^{(1)}_{u_1^{(1)},k^{(1)}_1}, \dots, \xx^{(d)}_{u_1^{(d)},k^{(d)}_1}) 
  \cdot \\ 
  &\phantom{=} \prod_{j=1}^{d} a^{(j)}(\fraku^{(j)},\kk^{(j)})
  \cdot 
  \left( \prod_{j=1}^{d} \chi^{(j)}_{u^{(j)}_{\# \fraku^{(j)}},k^{(j)}_
  {\# \fraku^{(j)}}} (\xx^{(j)}) \right),
    \end{aligned}
 \end{align}
 where 
 \begin{align*}
  c_{\bbeta}( \fraku^{(1)}, \dots, \fraku^{(d)}) = (-1)^{\| \bbeta \|_1 + d +
  \# \fraku^{(1)} + \cdots + \# \fraku^{(d)}}
 \end{align*}
and $ a^{(j)}(\fraku^{(j)}, \kk^{(j)}) = 1 $ if $ \# \fraku^{(j)} = 1 $ and 
\begin{align}\label{eq:combinedOperatorCoefficient}
a^{(j)}(\fraku^{(j)}, \kk^{(j)}) = \prod_{m = 1}^{\# \fraku^{(j)} - 1} 
\chi_{u^{(j)}_m,k^{(j)}_m} (\xx^{(j)}_{u^{(j)}_{m-1}, k^{(j)}_{m-1}}) 
\end{align}
for $ \#\fraku^{(j)} > 1 $ for every $ 1 \leq j \leq d $.
\end{definition}

\begin{remark}\label{rem:OnTPML}
\begin{enumerate}
 \item We can see that the target function is only evaluated in the sites $ \xx^{(j)}_{u_1^{(j)}, k_1^{(j)}} $,
 $ 1 \leq k_1^{(j)} \leq N_{u_1^{(j)}}^{(j)} $.
 These are the points in the coarsest set of sites $ X^{(j)}_{u_1^{(j)}} $ associated to each ordered set $ \fraku $.
 \item Similarly, we see that only the Lagrange functions associated to the finest set of sites for each ordered set $ \fraku $ are evaluated at the point $ \xx $.
\end{enumerate}
\end{remark}

Although this is a valid representation of the approximation operator, a naive implementation is, even in only two directions and a moderate number of levels in each direction, numerically too expensive.

\section{Modifications of the Tensor Product Multilevel Method}\label{sec:ModificationsOfTPML}

The aim is now to reduce the online cost of a point-evaluation of the tensor 
product multilevel approximation. We present two ways to achieve this goal. 
First, by moving expensive computations to an offline phase. And second, by 
deriving a completely new representation of the approximation operator, where
we change the perspective and sum of the points in the associated sparse grid
rather than the multi-indices.

\subsection{Towards more efficient point-evaluations}\label{sec:Efficient
Evaluations}

To keep this section brief we move the proofs into \cref{app:Proofs}. Most of 
the ideas are technical transformations involving the elements of the ordered
sets $ \fraku^{(j)} $.

We start with simplifying the representation of the tensor product of the combined
operators $ \cI^{(1)}_{\fraku^{(1)}} \otimes \cdots \otimes 
\cI^{(d)}_{\fraku^{(d)}} $. To do this, we take a closer look at each of the
direction-wise, combined operator $ \cI^{(j)}_{\fraku^{(j)}} $ and try to 
separate the intermediate terms in $ a^{(j)}(\fraku^{(j)}, \kk^{(j)}) $. These
can be pre-computed since they are independent of the target function and the 
evaluation point. This idea can be applied in each direction, hence we omit
the superscript $ (j) $ in the next lemma.

\begin{lemma}\label{lem:CombinedOperatorNewRepresentation}
Let $ \fraku \subseteq \{1, \dots, L \} $ be a fixed ordered set and 
 \begin{align*}
  \aalpha_{i,k} := \begin{cases} 
				  M_i^{-1} \ee_k, &\text{for interpolation,} \\
                  (M_i + \lambda_i I)^{-1} \ee_k, &\text{for penalized least-squares approximation}
                 \end{cases}
 \end{align*}
be the coefficient vector for $ \chi_{i,k} $.

Define $ P_{\fraku} \in  \R^{N_{u_1} \times N_{u_{\#\fraku}}} $ as
 \begin{align}\label{eq:ProductMatrix}
 P_{\fraku} := \prod_{m=2}^{\#\fraku} \sum_{k_m \leq N_{u_m}} \rr_{
 u_{m-1}}(\xx_{u_m,k_m}) \aalpha^{\transpose}_{u_m,k_m},
 \end{align}
with $ \rr_i $ as in \eqref{eq:definitionRVector}.

 Then we can express the combined local operators $ \cI_{\fraku} $ as 
 \begin{align*}
  \cI_{\fraku}f(\xx) = \sum_{k_1 \leq N_{u_1}} f(\xx_{u_1, k_1}) \aalpha
  ^{\transpose}_{u_1,k_1} \cdot P_{\fraku} \cdot \rr_{u_{\# \fraku}}(\xx),
  \quad \xx \in \Omega.
 \end{align*}
\end{lemma}

\begin{remark}\label{rem:FirstComputationalCost}
 \begin{enumerate}
 \item To compute the vectors $ (\aalpha_{i,k})_{1 \leq k \leq N_i} \in 
 \R^{N_i} $ we have to solve $ N_i $-many sparse linear systems, where the
 matrix stays the same only the right-hand side changes.
  \item We recall that $ \rr_i (\xx) $ is the evaluation vector of the rescaled 
  kernel on level $ i $. It has only constant many non-zero entries. This means that each of the rank-$ 1 $-matrices 
  \begin{align*}
   \rr_{u_{m-1}}(\xx_{u_m,k_m}) \aalpha^{\transpose}_{u_m,k_m}
  \end{align*}
  in \eqref{eq:ProductMatrix} is sparse. However, we cannot leverage this since for every $ k_m $ the sparsity pattern differs.
   \item All matrices $ P_{\fraku} $ can be computed in an offline phase and are 
 independent of the evaluation point $ \xx $ and the data 
 $ f(\xx_{u_1,k_1}) $. Additionally, all these computations are vector or matrix computations and, therefore, can be done efficiently on GPUs.
 \end{enumerate}
\end{remark}

We can now use \cref{lem:CombinedOperatorNewRepresentation} for every
direction $ 1 \leq j \leq d $ and obtain a new representation of the tensor 
product combined operator.

\begin{theorem}\label{thrm:NewRepresenationTensorProductCombinedOperator}
Use the notation and assumptions of \cref{def:TPMLGeneralRepresentation}. 
Additionally, assume that for every $ 1 \leq j \leq d $ and every ordered set
$ \fraku^{(j)} $ we use the notation of
\cref{lem:CombinedOperatorNewRepresentation}. Set 
\begin{align*}
\xx_{\uu_1,\kk_1} := (\xx_{u_1^{(1)},k_1^{(1)}}^{(1)}, \dots, \xx_{u_1^
  {(d)},k_1^{(d)}}^{(d)}).
\end{align*}
 
Then the tensor product combined operator can be expressed as
 \begin{align}\label{eq:NewRepresenationTensorProductCombinedOperator}
 \begin{aligned}
 & \hspace{-3em}\left(\cI_{\fraku^{(1)}}^{(1)} \otimes \cdots \otimes 
  \cI_{\fraku^{(d)}}^{(d)} \right) (\ff)(\xx) =  \\
&= \sum_{k_1^{(1)} \leq N_{u_1^{(1)}}^{(1)}} \cdots \sum_{k_1^{(d)} \leq 
  N_{u_1^{(d)}}^{(d)}} \ff(\xx_{\uu_1,\kk_1}) \cdot  \\
&\phantom{5em} \cdot \prod_{j=1}^{d} ( \aalpha^{(j)}_{u_1^{(j)},
  k_1^{(j)}} 
  )^{\transpose} \cdot P^{(j)}_{\fraku^{(j)}} \cdot \rr^{(j)}_{u^{(j)}
  _{\# \fraku^{(j)}}} ( \xx^{(j)}). 
\end{aligned}
\end{align}
\end{theorem}

The new representation of $ \left(\cI_{\fraku^{(1)}}^{(1)} \otimes \cdots
\otimes \cI_{\fraku^{(d)}}^{(d)} \right) (\ff)(\xx) $ immediately yields a simplified intermediate representation of the multilevel tensor product approximation.

\begin{corollary}\label{cor:IntermediateNewTPML}
With the notation and assumptions of 
\cref{thrm:NewRepresenationTensorProductCombinedOperator} the tensor product
multilevel approximant $ \cA_{\cI_{\oomega}(\ell,d)}(\ff) $ can, for every 
$ \xx \in \Omega $, be expressed as
\begin{align}\label{eq:TPMLWithFirstSimplifications}
\begin{aligned}
 &\cA_{\cI_{\oomega}(\ell,d)}(\ff)(\xx) = \sum_{\llambda \in \cJ_{\oomega}
 (\ell,d)} \sum_{\substack{\bbeta \in \{0,1\}^d \\ \llambda + \bbeta \in \cI_
 {\oomega}(\ell,d)}}(-1)^{\|\bbeta\|_1 + d} \\
 & \sum_{u^{(1)}_1 \in \{1, \dots, \lambda_1 \}} \cdots \sum_{u^{(d)}_1 \in
 \{1, \dots, \lambda_d \}} \sum_{k_1^{(1)} \leq N^{(1)}_{u_1^{(1)}}} \cdots
 \sum_{k_1^{(d)} \leq N^{(d)}_{u_1^{(d)}}} \ff(\xx_{u_1, k_1}) \cdot \\
 &\prod_{j=1}^{d} \left(\aalpha^{(j)}_{u_1^{(j)},k_1^{(j)}}\right)^{\transpose}
 \sum_{\substack{ \widetilde{\fraku}^{(j)} \subseteq \{1,\dots,\lambda_j\} \\ 
 \widetilde{u}^{(j)}_1 = u_1^{(j)}}} (-1)^{\# \widetilde{
 \fraku}^{(j)}}  P^{(j)}_{\widetilde{\fraku}^{(j)}} \cdot 
  \rr^{(j)}_{\widetilde{u}^{(j)}_{\# \widetilde{\fraku}^{(j)}}} ( \xx^{(j)}).
\end{aligned}
\end{align}
\end{corollary}

\begin{remark}
 The ordered sets $ \widetilde{\fraku}^{(j)} $ in the inner most sum in 
 \eqref{eq:TPMLWithFirstSimplifications} are those ordered sets 
 $ \fraku^{(j)} $ of \eqref{eq:TPMLGeneralRepresentation} whose first element
 is the, for the inner sum fixed, $ u_1^{(j)} $.
\end{remark}

The final step is now to simplify the inner most sum in
\eqref{eq:TPMLWithFirstSimplifications}. This will be achieved by also fixing
the last element of the ordered sets $ \widetilde{\fraku} $. This allows us to pre-compute most of the intermediate sums over the matrices $ P^{(j)}_
 {\widetilde{\fraku}^{(j)}} $.

\begin{corollary}\label{cor:InnermostSum}
For $ 1 \leq m \leq p \leq \lambda_j $, $ 1 \leq j \leq d $, set 
\begin{align}\label{eq:SumMatrix}
S^{(j)}_{m,p} := \sum_{\substack{\overline{\fraku}^{(j)} 
\subseteq \{1, \dots, \lambda_j\} \\ \overline{u}^{(j)}_1 = m \\
\overline{u}^{(j)}_d = p}} (-1)^{\# \overline{\fraku}^{(j)}} P^{(j)}_
{\overline{\fraku}^{(j)}} \quad \in \R^{N^{(j)}_m \times N^{(j)}_p}.
\end{align}
Then we have
\begin{align}\label{eq:InnermostSum}
 &\sum_{\substack{ \widetilde{\fraku}^{(j)} \subseteq \{1,\dots,\lambda_j\} \\ 
 \widetilde{\fraku}^{(j)}_1 = u_1^{(j)}}} (-1)^{\# \widetilde{
 \fraku}^{(j)}}  P^{(j)}_{\widetilde{\fraku}^{(j)}} \cdot \rr^{(j)}
 _{\widetilde{u}^{(j)}_{\# \widetilde{\fraku}^{(j)}}} ( \xx^{(j)}) = 
\sum_{p = u_1^{(j)}}^{\lambda_j} S^{(j)}_{u_1^{(j)}, p} \cdot \rr^{(j)}_{p}
  (\xx^{(j)}).
\end{align}
\end{corollary}

\begin{remark}
\begin{enumerate}
    \item The matrix $ S^{(j)}_{m,p} $ in \eqref{eq:SumMatrix} considers all matrices $ P_{\fraku} $ whose ordered set $\overline{\fraku}^{(j)}$ has $ m $ as first and $ p $ as last entry.
    \item We can express the sum on the right-hand side of \eqref{eq:InnermostSum} as
    \begin{align*}
        \left(S^{(j)}_{u^{(j)}_1,u^{(j)}_1} \ S^{(j)}_{u^{(j)}_1, u^{(j)}_1 + 1} \cdots S^{(j)}_{u^{(j)}_1, \lambda_j} \right) \begin{pmatrix} \rr^{(j)}_{u^{(j)}_1}(\xx^{(j)}) \\ \rr^{(j)}_{u^{(j)}_1 + 1}(\xx^{(j)}) \\ \vdots \\ \rr^{(j)}_{u^{(j)}_{\lambda_j}}(\xx^{(j)})
         \end{pmatrix}.
    \end{align*}
    \item We collect the matrices $ S^{(j)}_{m,p} $ for all $ 1 \leq m \leq p \leq \lambda_{j,max} $ in the system
    \begin{align}\label{eq:LargeSMatrix}
        \SS^{(j)} :=
        \begin{pmatrix} 
        S^{(j)}_{1,1} & S^{(j)}_{1,2} & \hdots & S^{(j)}_{1, \lambda_{j,max}} \\ 
                      & S^{(j)}_{2,2} & \hdots & S^{(j)}_{2,\lambda_{j,max}}  \\
                      &               & \ddots & \vdots                       \\
                      &               &        & S^{(j)}_{\lambda_{j,max},\lambda_{j,max}}
        \end{pmatrix}.                  
    \end{align}
\end{enumerate}
\end{remark}

Finally, we can interpret the term in the product in \eqref{eq:TPMLWithFirstSimplifications},
\begin{align*}
    \left(\aalpha^{(j)}_
 {u_1^{(j)},k_1^{(j)}}\right)^{\transpose}
 \sum_{p = u_1^{(j)}}^{\lambda_j}S^{(j)}_{u_1^{(j)}, p} \cdot \rr^{(j)}_{p}
  (\xx^{(j)})
\end{align*}
as a product of block-matrices. Using the notation 
\begin{align*}
    \AA^{(j)}_i = \begin{pmatrix} \aalpha^{(j)}_{i,1} \\ \vdots \\ \aalpha^{(j)}_{i,N_i^{(j)}} \end{pmatrix},
\end{align*}
we define the matrix $ \XXi^{(j)}(\xx^{(j)}) $ by
\begin{align}\label{eq:defXXi}
    &\XXi^{(j)}(\xx^{(j)}) := 
    \begin{pmatrix}
        \xxi^{(j)}_{1,1}(\xx^{(j)}) & \cdots & \xxi^{(j)}_{1, \lambda_{j,max}}(\xx^{(j)} \\
        & \ddots & \vdots \\
        0 & & \xxi^{(j)}_{\lambda_{j,max}, \lambda_{j,max}}(\xx^{(j)})
    \end{pmatrix} \nonumber \\
    &= \begin{pmatrix} 
    A^{(j)}_i &        & 0 \\
              & \ddots &   \\
    0         &        & A^{(j)}_{\lambda_{j,max}}
    \end{pmatrix}
    \cdot \SS^{(j)}
    \cdot
    \begin{pmatrix}
        \rr^{(j)}_{1} (\xx^{(j)})&         & 0 \\
                      & \ddots  &   \\
                    0 &         & \rr^{(j)}_{\lambda_{j,max}}(\xx^{(j)})
    \end{pmatrix}.
\end{align}

Clearly, the matrix $ \XXi^{(j)}(\xx^{(j)}) $ depends on the $ j $-th component of the evaluation point $ \xx $, by the third factor. We recall that every $ \rr^{(j)}_i(\xx^{(j)}) $ is a vector with only a constant number of nonzero entries. So, in application, one has to decide whether to compute the matrix as in \eqref{eq:defXXi}, by using the sparseness of the third matrix, for every evaluation point, or just compute the product of the first two dense matrices and keep the independence of the evaluation point.

Altogether, the tensor product multilevel interpolation or penalized 
least-squares approximation has the following form. This is the main result of 
this section and is the representation that allows the most pre-computation.

\begin{theorem}\label{thrm:FinalRepresentation}
With the notation and assumptions made throughout, we can represent the tensor product multilevel
approximation to $ \ff $ as
\begin{align}\label{eq:NewRepresentationTPMLOperator}
 &\cA_{\cI_{\oomega}(\ell,d)}(\ff)(\xx) = \nonumber \\
 &\sum_{\llambda \in \cJ_{\oomega}
 (\ell,d)} \sum_{\substack{\bbeta \in \{0,1\}^d \\ \llambda + \bbeta \in \cI_
 {\oomega}(\ell,d)}}(-1)^{\|\bbeta\|_1 + d} \nonumber \\
 &\phantom{\sum} \sum_{u^{(1)}_1 \in \{1, \dots, \lambda_1 \}} \cdots \sum_
 {u^{(d)}_1 \in \{1, \dots, \lambda_d \}} \sum_{k_1^{(1)} \leq N^{(1)}_
 {u_1^{(1)}}} \cdots \sum_{k_1^{(d)} \leq N^{(d)}_{u_1^{(d)}}} 
 \ff(\xx_{u_1, k_1}) \cdot \nonumber\\
 &\phantom{\sum \sum \sum}\cdot \prod_{j=1}^{d} \sum_{m = u_1^{(j)}}^{\lambda_{j}} \left( \xxi^{(j)}_{u_1^{(j)},m} (\xx^{(j)}) \right)_{k_1^{(j)}} 
\end{align}
for every $ \xx \in \Omega^{(1)} \times \cdots \times \Omega^{(d)} $.
\end{theorem}

\subsection{A Nodal Representation}\label{subsec:NodalRepresentation}

In the previous section, we found a representation of the operator \eqref{eq:TPMLGeneralRepresentation} by smartly finding ways to pre-compute certain terms. However, it still has several downsides. For one, we need the same value $ \ff(\xx_{u_1,k_1}) $ multiple times, which can be alleviated by storing all those values in an easy-to-access data structure. But more importantly, we have no advantage, if the direction-wise sites are nested, i.e., if 
\begin{align*}
    X^{(j)}_1 \subset X^{(j)}_2 \subset \cdots \subset X^{(j)}_{\lambda_{j,max}}.
\end{align*}
It turns out that in this case, we can find another representation of the operator in \eqref{eq:TPMLGeneralRepresentation}. This new representation has the advantage that the outer sum is over the points in the associated sparse grid. 

To do this, we need to introduce some more notation. In a nested family of sets of sites, for each point $ \xx \in X^{(j)}_{\lambda_{j,max}} $, there is a unique level $ u(\xx) \in \{1, \dots, \lambda_{j,max} \} $ where it first occurs. This point has position $ k_j(\xx,i_j) $ in the sets $ X^{(j)}_{i_j} $, $ u(\xx) \leq i_j \leq \lambda_{j,max} $.

With this notation, we can find a new representation of the direction-wise multilevel operators.

\begin{theorem}\label{thrm:NodalRepresentationMultilevel}
With the notation and assumptions made throughout this paper, we define
\begin{align*}
    A_{\{u(\xx), \dots, L\}} \chi_{L,k(\xx,L)} := \sum_{\emptyset \neq \fraku \subseteq\{ u(\xx), \dots, L \}} (-1)^{\#\fraku + 1} \cI_{\fraku} \chi_{L,k(\xx,L)}.
\end{align*}
Then the multilevel operator introduced in \cref{thrm:RepresentationMultilevelOperator} can be expressed as
\begin{align*}
    A_L(f) := \sum_{\yy \in X_L} f(\yy) A_{\{u(\yy), \dots, L\}} \chi_{L,k(\yy,L)}.
\end{align*}
\end{theorem}

Using this representation of the direction-wise multilevel operators in the combination technique and easy manipulations of the sums yield the \emph{nodal representation} of the tensor product multilevel operator.

\begin{theorem}\label{thrm:NodalRepresentationTPML}
    With the notation and assumption made throughout this paper, we can represent the tensor product multilevel operator $ \cA_{\cI_{\oomega}(\ell,d)} $, defined in \cref{def:TPMLGeneralRepresentation}, as
    \begin{align*}
        \cA_{\cI_{\oomega}(\ell,d)}(\ff) &= \sum_{\xx \in \cH_{\cI_{\oomega}(\ell,d)}} \ff(\xx) \cdot \\
        &\phantom{=} \cdot\sum_{\substack{\llambda \in \cJ_{\oomega}(\ell,d) \\ \llambda \geq \uu(\xx)}} \sum_{\substack{\bbeta \in \{0,1 \}^d \\ \llambda + \bbeta \in \cI_{\oomega}(\ell,d)}} (-1)^{\| \bbeta\|_1 + d} \cdot \\
        &\phantom{\sum \sum} \cdot \bigotimes_{j=1}^{d} A^{(j)}_{\{u_j(\xx^{(j)}), \dots, \lambda_j \}} \chi^{(j)}_{\lambda_j,k_j(\xx^{(j)},\lambda_j)}.
    \end{align*}
\end{theorem}

\begin{remark}
\begin{enumerate}
    \item As a direct consequence, the nodal representation of the tensor product multilevel method in \cref{thrm:NodalRepresentationTPML} is an ideal candidate for GPU parallelization: For every point in the sparse grid, we do more or less the same computations.
    \item Again, we can pre-compute all the direction-wise multilevel approximations $ A^{(j)}_{\{ u_j(\xx^{(j)}), \dots, \lambda_j \}} \chi^{(j)}_{\lambda_j,k_j(\xx^{(j)},\lambda_j)} $ and keep them in memory. For a single evaluation of the approximation at a point $ \yy = (\yy^{(1)}, \dots, \yy^{(d)}) \in \Omega^{\otimes} $ we need to evaluate the stored multilevel approximations once and store them in an appropriate data structure such that getting the value $ A^{(j)}_{\{ u_j(\xx^{(j)}), \dots, \lambda_j \}} \chi^{(j)}_{\lambda_j,k_j(\xx^{(j)},\lambda_j} (\yy^{(j)})$ is a simple lookup.
\end{enumerate}
\end{remark}

\section{Numerical Examples}\label{sec:NumericalExamples}

We will now showcase two example applications for the tensor product multilevel method. 
In the first test case, the method is applied to interpolate simulation data from a shallow water simulation in space and time. There, the tensor product multilevel method is setup with a two-dimensional domain (the spatial domain of the simulation) in the first direction and a one-dimensional domain (the time) in the second direction.

The second example interpolates a seven-dimensional function on a sparse grid. The interpolated function in this case is the shear force of a cantilever beam, and the tensor product multilevel method is set up with completely isotropic directions.

\subsection{Tidal Flow at the Bight of Abaco}\label{subsex:TidalFlow}
In the first example, we apply the tensor product multilevel method to interpolate simulation data from a tidal flow simulation at the Bight of Abaco. 
The tidal flow is simulated by the two-dimensional shallow water equations, given by
\begin{equation}
\label{eq:swe}
\begin{aligned}
&\partial_t \xi+\nabla \cdot \boldsymbol{q}=0,\\
&\partial_t \boldsymbol{q}
+\nabla \cdot \left(\boldsymbol{qq}^T /H \right)
+\tau_{\textrm{bf}} \boldsymbol{q}
+\left( \begin{smallmatrix} 0 & -f_c\\ f_c & 0 \end{smallmatrix} \right) \boldsymbol{q}
+gH\nabla\xi
=\boldsymbol{F},
\end{aligned}
\end{equation}
where $\xi$ is the water height above some median sea level, $\boldsymbol{q} = (U, V)^T$ denotes the depth-integrated horizontal velocities, and $H = \xi + h_b$ is the total water depth relative to the bathymetry $h_b$.
Other variables in the model are the Coriolis coefficient $f_c$, the bottom friction coefficient $\tau_bf$, the gravitational acceleration $g$, and external forces $\boldsymbol{F}$.

\begin{figure}
    \centering
    \includegraphics[width=0.5\linewidth]{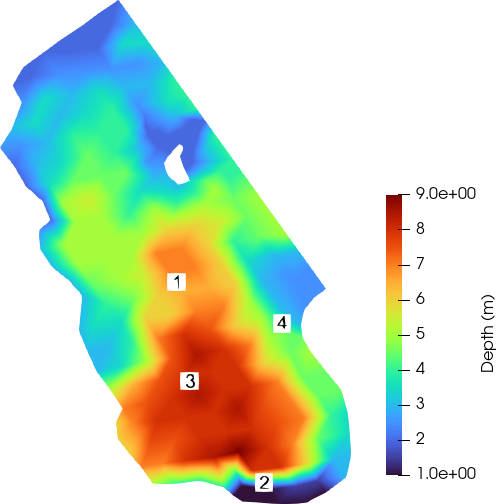}
    \caption{Bathymetry for the Bight of Abaco. Horizontal extends of the domain are 70 km in $x$ direction and 100 km in $y$ direction. Four virtual recording stations report the variables $\xi$, $U$ and $V$ at the positions marked by 1 -- 4.}
    \label{fig:bahamas-bathymetry}
\end{figure}

Equation \eqref{eq:swe} is discretized with a discontinuous Galerkin method on unstructured triangular meshes. Details of the discretization and implementation can be found in \cite{buettner:ShallowWaterDG}. 
Figure \ref{fig:bahamas-bathymetry} shows the bathymetry of the Bight of Abaco, the exact boundary conditions and parameters for this domain can be found in \cite[Section 2.3]{buettner:ShallowWaterDG}. 
We run a tidal simulation with this bathymetry on a mesh with approximately 55,000 grid points for a simulation period of three days. After an initial ramp-up of two days for the boundary conditions, the variables $\xi$, $U$, and $V$ are written for all grid points at intervals of five minutes during the third simulation day. This data then serves as input to the multilevel method.

The tensor product multilevel method is now set up with two directions. The first direction represents the two-dimensional spatial coordinate within the domain, and the second direction represents the time. 
We generate a set of nested points for the spatial direction by iteratively thinning the mesh points, so that the separation distances $q_i$ satisfy $q_{i+1} \approx \frac{1}{2} q_i$ for levels $i = 1, \dots, 6$.
On the coarsest level, this results in merely 35 points in the whole domain, whereas the finest level contains approximately 38,000 points.
We use the compactly supported kernel $\phi_{3,1}(r) = (1-r)^4_+(4r+1)$
for the RBF interpolation, and the support radius $\varepsilon_i$ is coupled to the separation distance at each level by $\varepsilon_i = 6 q_i$.

For the time direction, we used the initial output frequency of five minutes on the finest level and double the time step for each coarser level. As a basis function for interpolation, we use the kernel $\phi_{1,1}(r) = (1-r)^3_+(3r+1)$ with the support radius set to $\varepsilon_i = 6 \Delta t_i$. Table \ref{tab:SweMultilevelSets} summarizes the number of points, separation distance, and time step $\Delta t$ for all levels.

\begin{table}[t]
    \centering
    \caption{Number of points, separation distance and time step for each level for the Bight of Abaco simulation.}
    \begin{tabular}{l rrrrrr}
    \toprule
    Level & 1 & 2 & 3 & 4 & 5 & 6 \\ 
    Number of points    &     35 &     123 &    476 &   1871 &  8228 & 37953 \\
    Separation distance (m) & 9133.3 & 4521.2 & 2254.9 & 1128.1 & 563.7 & 281.9 \\
    Points in time direction &      9 &      18 &     36 &     72 &   144 &   288 \\
    $\Delta t$ (minutes)&    160 &      80 &     40 &     20 &    10 &     5 \\ 
    \bottomrule
    \end{tabular}
    \label{tab:SweMultilevelSets}
\end{table}

\begin{figure}
    \centering
    \includegraphics[height=8cm]{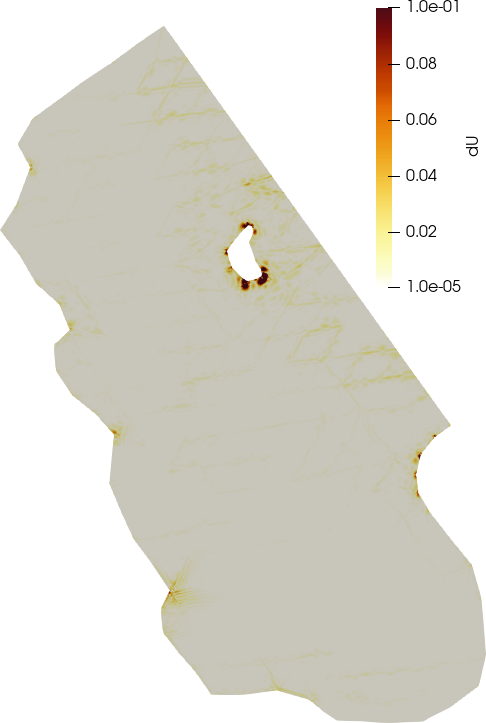}
    \hfill
    \includegraphics[height=8cm]{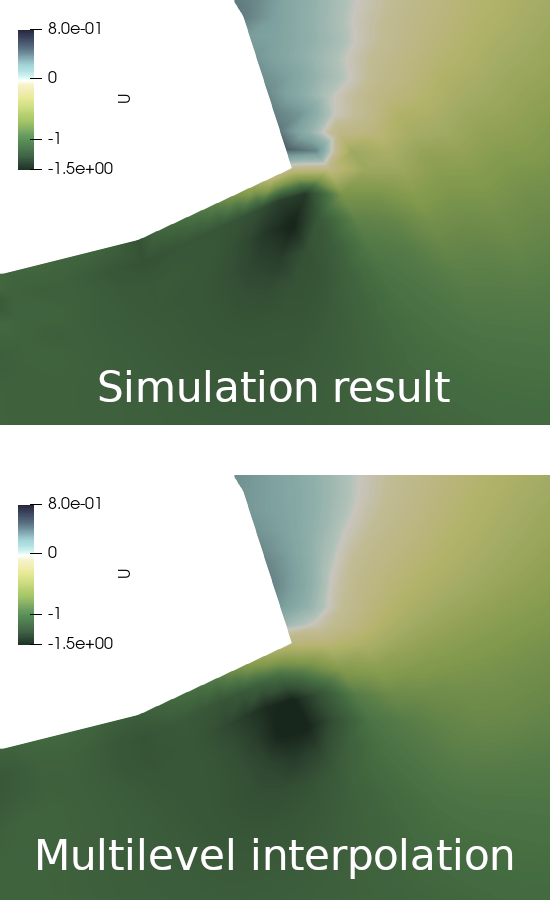}
    \caption{Absolute pointwise difference for $U$ on the whole mesh (left), simulation and interpolation result for $U$ near the island (right).}
    \label{fig:SweErrorU}
\end{figure}

We use the mesh with 55,000 points, which generated the original data, to evaluate the multilevel method.
By evaluating the interpolant on these points and comparing the result point-wise with the simulation data, we are able to quantify the error of the multilevel method.
Figure \ref{fig:SweErrorU} (left) shows the difference between the interpolant and the simulation result for the variable $U$. 
The largest differences appear near the boundary of the domain, especially near the island. Zooming in on this particular region reveals that the interpolated data looks much smoother than the simulation data, as even element boundaries are visible there (Figure~\ref{fig:SweErrorU}, top right). Apart from this difference, a statistical analysis shows that for 95~\% of the mesh points the difference lies below $10^{-2}$.

\begin{figure}
    \centering
    \begin{tikzpicture}
        \begin{axis}[ymode=log, xlabel={Number of levels}, ylabel={Average pointwise error}]
            \addplot coordinates{(1, 0.024)(2, 0.00513)(3, 0.00148)(4, 0.000303)(5, 0.0000638)(6, 0.0000207)};
            \addplot coordinates{(1, 0.214)(2, 0.049)(3, 0.0198)(4, 0.00599)(5, 0.00205)(6, 0.0007982)};
            \addplot coordinates{(1, 0.151)(2, 0.0406)(3, 0.0154)(4, 0.00467)(5, 0.00167)(6, 0.0007875)};
            \legend{$\xi$, $U$, $V$}
        \end{axis}
    \end{tikzpicture}
    \caption{Average pointwise error for the Bahamas example with different number of levels and variables.}
    \label{fig:SweErrorAverage}
\end{figure}
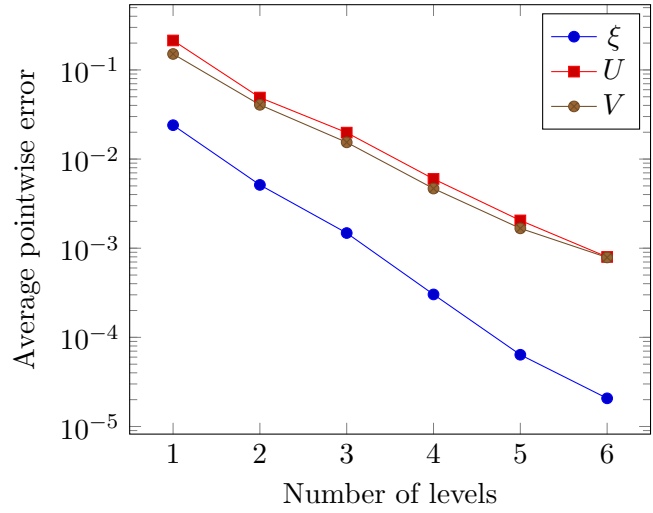

Figure~\ref{fig:SweErrorAverage} shows the average pointwise error for different number of levels and all three variables. 
The water height $\xi$ is approximated better than the velocities and has a higher convergence order: For $\xi$, the computed order of convergence is approximately 2.0, whereas the convergence order of $U$ and $V$ is 1.6 and 1.5 respectively.

\begin{figure}
    \centering
    \includegraphics[width=0.5\linewidth]{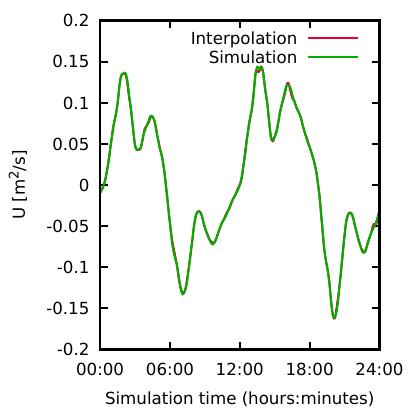}%
    \includegraphics[width=0.5\linewidth]{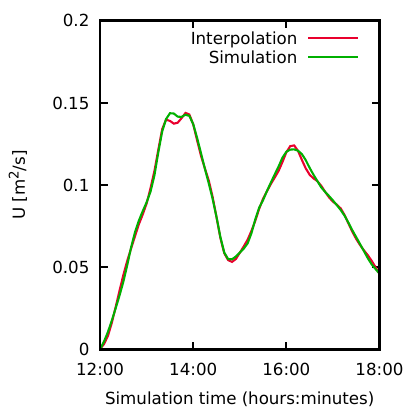}
    \caption{Depth-integrated velocity in $x$ direction at recording station 2.}
    \label{fig:BahamasStation2}
\end{figure}

If we focus on a single recording station (\cref{fig:BahamasStation2} left), we see that the velocity profile from the interpolant matches the simulation data. However, there are also some slight deviations visible (\cref{fig:BahamasStation2} right). Because the station locations are within the elements, and not part of the interpolation sites, a slight deviation is not surprising.

\subsection{Cantilever beam}
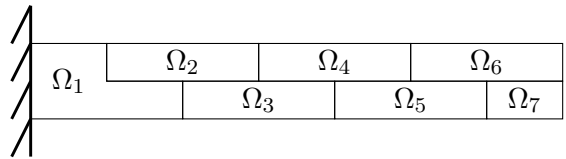
\begin{figure}
	\centering
	\begin{tikzpicture}[y=-1cm]
		\draw (0, 0) -- (1, 0) -- (1, 0.5) -- (2, 0.5) -- (2, 1) -- (0, 1) -- cycle;
		\draw (1, 0) -- (3, 0) -- (3, 0.5) -- (1, 0.5) -- cycle;
		\draw (2, 0.5) -- (4, 0.5) -- (4, 1) -- (2, 1) -- cycle;
		\draw (3, 0) -- (5, 0) -- (5, 0.5) -- (3, 0.5) -- cycle;
		\draw (4, 0.5) -- (6, 0.5) -- (6, 1) -- (4, 1) -- cycle;
		\draw (5, 0) -- (7, 0) -- (7, 0.5) -- (5, 0.5) -- cycle;
		\draw (6, 0.5) -- (7, 0.5) -- (7, 1) -- (6, 1) -- cycle;
		\node at (0.5, 0.5) {$\Omega_1$};
		\node at (2, 0.25) {$\Omega_2$};
		\node at (3, 0.75) {$\Omega_3$};
		\node at (4, 0.25) {$\Omega_4$};
		\node at (5, 0.75) {$\Omega_5$};
		\node at (6, 0.25) {$\Omega_6$};
		\node at (6.5, 0.75) {$\Omega_7$};
		\draw[line width=1.1] (0, -0.5) -- (0, 1.5);
		\draw[line width=1.1] (0, -0.5) -- (-0.25, 0);
		\draw[line width=1.1] (0, 0) -- (-0.25, 0.5);
		\draw[line width=1.1] (0, 0.5) -- (-0.25, 1);
		\draw[line width=1.1] (0, 1) -- (-0.25, 1.5);
	\end{tikzpicture}
	\caption{Cantilever beam from \cite[Section 3.4]{migliorati:ApproximationQOI}.}
	\label{fig:cantilever}
\end{figure}

As a second test example, we consider a cantilever beam, as described in \cite{migliorati:ApproximationQOI} and depicted in Figure~\ref{fig:cantilever}. The displacement field $\uu$ satisfies the Navier-Lam\'e equation,
\begin{align*}
	(\lambda(\xx, \yy) + \mu(\xx, \yy))\nabla(\nabla \cdot \uu(\xx, \yy)) + \mu(\xx, \yy)\nabla^2 \uu(\xx, \yy) &= -\ff(\xx, \yy), && \xx \in \Omega \\
	\sigma(\uu(\xx, \yy)) \cdot \nn &= 0, && \xx \in \partial\Omega \setminus \partial\Omega_{wall} \\
	\uu(\xx, \yy) &= \boldsymbol{0}, && \xx \in \partial\Omega_{wall},
\end{align*}
where $\mu$ and $\lambda$ are the Lam\'e constants given by
\begin{equation*}
	\mu(\xx, \yy) = \frac{E(\xx, \yy)}{2(1 + \nu)} \quad \text{ and } \quad \lambda(\xx, \yy) = \frac{\nu E(\xx, \yy)}{(1 + \nu)(1 - 2\nu)},
\end{equation*}
$\nu = 0.28$ is Poisson's ratio, $E$ is Young's modulus and $\sigma$ is the Cauchy stress tensor computed by
\begin{equation*}
	\sigma(\uu(\xx, \yy)) = \lambda\left(\nabla \cdot \uu(\xx, \yy)\right) I + \mu \left( \nabla \uu(\xx, \yy) + (\nabla \uu(\xx, \yy))^T\right).
\end{equation*}
As shown in the figure, the beam is separated into seven non-overlapping subdomains $\Omega_i$, $i = 1,\dots, 7$, each having a different value for Young's modulus:
\begin{equation*}
	E(\xx, \yy) = \exp(7 + y_i) \text{ if } \xx \in \Omega_i,~i = 1, \dots, 7.
\end{equation*}

\begin{figure}
    \centering
    \begin{tikzpicture}
        \begin{axis}[ymode=log, xlabel={Number of levels}, 
            title={Maximum relative error},width=0.5\linewidth]
            \addplot coordinates{(1, 0.128)(2, 0.0432)(3, 0.00946)(4,0.00208)(5,0.000561)(6,0.000117)};
            \addplot coordinates{(1, 0.125)(2, 0.0319)(3, 0.00681)(4,0.00133)(5,0.000338)(6,6.52e-5)};
            \legend{$\phi_{1,1}$, $\phi_{1,2}$}
        \end{axis}
    \end{tikzpicture}\begin{tikzpicture}
        \begin{axis}[ymode=log, xlabel={Number of levels}, 
            title={Maximum absolute error},width=0.5\linewidth]
            \addplot coordinates{(1, 2.954)(2,1.136)(3,0.283)(4,0.0557)(5,0.0148)(6,0.00336)};
            \addplot coordinates{(1, 3.11)(2,0.795)(3,0.169)(4,0.033)(5,0.00841)(6,0.00162)};
            \legend{$\phi_{1,1}$, $\phi_{1,2}$};
        \end{axis}
    \end{tikzpicture}
    \caption{Maximum relative and absolute error for the beam model with different number of levels.}
    \label{fig:BeamRelativeError}
\end{figure}
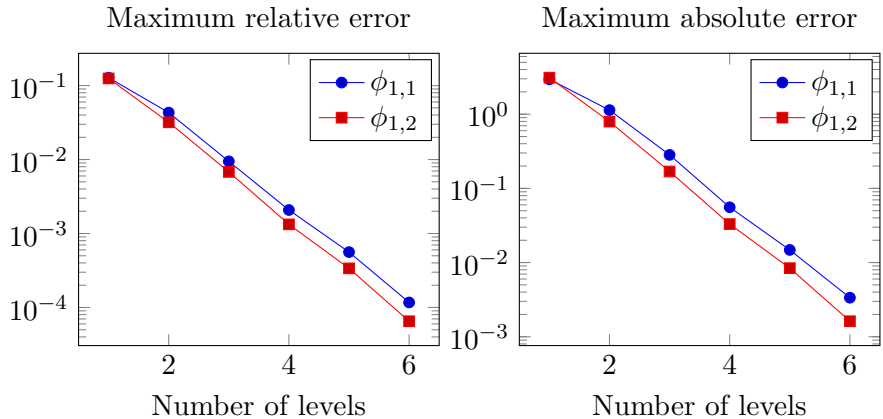

We set, following the description in \cite{migliorati:ApproximationQOI}, $f \equiv 1$ and compute the total shear force in $y$-direction by
\begin{equation}
    q(\yy) = \int_\Omega \sigma_{12}(\uu(\xx, \yy) \mathop{d\xx},
\end{equation}
and use the tensor product multilevel method to approximate this function. 
We decompose the domain $[-1, 1]^7$ into one-dimensional intervals, and generate a sequence of nested equidistant points in each direction, where level $i$ contains $2^i + 1$ points. These points are then used in combination with the compactly supported basis functions
\begin{align*}
    \phi_{1,1}(r) &= (1-r)^3_+(3r+1) \\
    \phi_{1,2}(r) &= (1-r)^5_+(8r^2+5r+1)
\end{align*}
for the tensor product multilevel method. The support radius is set to four times the separation distance. 
Figure \ref{fig:BeamRelativeError} shows the maximum relative and absolute error for the beam model for both kernels.

\bigskip\noindent
{\bf Acknowledgement.}
The work of Holger Wendland was partially funded by the Deutsche Forschungsgemeinschaft (DFG, German Research Foundation) - Project Number 452806809.

We thank Vadym Aizinger for kindly providing the data used in \cref{subsex:TidalFlow}.

\newpage
\appendix

\section{Proofs of the statements of \cref{sec:Efficient
Evaluations}}\label{app:Proofs}

\begin{proof} [Proof of \cref{lem:CombinedOperatorNewRepresentation}]
 We recall the representation of $ \cI_{\fraku} $ in 
 \eqref{eq:defCombinedOperator},
 \begin{align*}
  \cI_{\fraku} = \sum_{\kk \leq \NN_{\fraku}} f(\xx_{u_1,k_1}) a(\fraku,\kk) 
\chi_{\fraku_{\# \fraku},k_{\# \fraku}}.
 \end{align*}
 Inserting the representations for the Lagrange functions \eqref{eq:LagrangeFunctionVectorVectorProduct}, yields
 \begin{align*}
  &\cI_{\fraku}(f)(\xx) =\\
  &= \sum_{\kk \leq \NN_{\fraku}} f(\xx_{u_1, k_1})
  \left[ \prod_{m=1}^{\# \fraku - 1} \aalpha^{\transpose}_{u_m,k_m} \rr_
  {u_m}(\xx_{u_{m+1},k_{m+1}}) \right] \aalpha^{\transpose}_{u_{\#\fraku},
  k_{\#\fraku}} \rr_{u_{\#\fraku}} (\xx).
 \end{align*}
 We see that in the product in the brackets, the index $ u_1,k_1 $ only 
 appears in $ \aalpha_{u_1,k_1} $. Similarly, in the product, there is no
 $ \aalpha^{\transpose}_{u_{\#\fraku},k_{\#\fraku}} $. This vector is written 
 explicitly after the bracket. Hence, we can split $ \aalpha^{\transpose}
 _{u_1,k_1} $ from the product and add $ \aalpha^{\transpose}
 _{u_{\#\fraku},k_{\#\fraku}} $ to it. This yields,
 after an index shift in the product term,
 \begin{align*}
  &\cI_{\fraku}(f)(\xx) =\\
  &= \sum_{\kk \leq \NN_{\fraku}} f(\xx_{u_1, k_1}) 
  \aalpha^{\transpose}_{u_1,k_1}
  \left[ \prod_{m=2}^{\# \fraku} \rr_{u_{m-1}}(\xx_{u_m,k_m}) 
  \aalpha^{\transpose}_{u_m,k_m}\right] \rr_{u_{\#\fraku}} (\xx).
 \end{align*}
 Next, we split the multiple sums $ \sum_{\kk \leq \NN_{\fraku}} $ in its 
 single components. Only the terms $ f(\xx_{u_1, k_1}) $ and $ 
  \aalpha^{\transpose}_{u_1,k_1} $ depend on $ k_1 $. We pull the other sums
  into the product. We have
  \begin{align*}
&\cI_{\fraku}(f)(\xx) = 
\sum_{k_1 \leq N_1} f(\xx_{u_1, k_1}) \aalpha^{\transpose}_{u_1,k_1}\cdot \\
&\phantom{\cI_{\fraku}(f)(\xx) = \sum} \cdot
\left[ \prod_{m=2}^{\# \fraku} \sum_{k_m \leq N_m} \rr_{u_{m-1}}
  (\xx_{u_m,k_m})
   \aalpha^{\transpose}_{u_m,k_m}\right] 
  \rr_{u_{\#\fraku}} (\xx).
  \end{align*}
 This is the claim if we define $ P_{\fraku} $ as in \eqref{eq:ProductMatrix}.
\end{proof}

\begin{proof}[Proof of \cref{thrm:NewRepresenationTensorProductCombinedOperator}]
The statement in \cref{thrm:NewRepresenationTensorProductCombinedOperator} follows easily from the observation that 
\begin{align*}
   \left( \chi^{(1)}_{i_1} \otimes \cdots \otimes \chi^{(d)}_{i_d} \right) (\xx) = \prod_{j=1}^{d} \chi^{(j)}_{i_j} (\xx^{(j)}) 
\end{align*}
 and then applying \cref{lem:CombinedOperatorNewRepresentation} in every direction independently.
\end{proof}

\begin{proof}[Proof of \cref{cor:IntermediateNewTPML}]
As already outlined in the text, this is an application of \cref{thrm:NewRepresenationTensorProductCombinedOperator}, together with using the observation that we only use evaluations of $ \ff $ in points on level $ u_1^{(j)} $. The claim is then obtained by clever grouping of the different summands.
\end{proof}

The proofs of the remaining \cref{cor:InnermostSum} and \cref{thrm:FinalRepresentation} were already outlined in the text.

\section{Proofs of the statements of \cref{subsec:NodalRepresentation}}\label{app:Proofs2}

\begin{proof}[Proof of \cref{thrm:NodalRepresentationMultilevel}]
Following \cref{thrm:RepresentationMultilevelOperator}, the multilevel operator $A_L$ has the representation
\begin{align*}
A_L(f) = \sum_{\substack{\fraku \subseteq \{1, \dots, L \} \\ 1 \leq \# \fraku \leq L}} (-1)^{\# \fraku + 1} \sum_{\kk \leq \NN_{\fraku}} a(\fraku, \kk) f(\xx_{u_1,k_1}) \chi_{u_{\# \fraku}, k_{\# \fraku}}.
\end{align*}
Again, we split off the first elements of the ordered sets $ \fraku $ and $ \kk $. yielding the decompositions $ \fraku = \{ u_1 \} \cup \widetilde{\fraku} $ and $ \kk = (k_1, \widetilde{\kk} ) $. Here, $ \widetilde{\fraku} \subseteq \{ u_1 + 1, \dots, L \} $ and $ \widetilde{\kk} \leq \NN_{\widetilde{\fraku}} $ and $ \widetilde{\fraku} = \emptyset $ is allowed and means that $ \fraku = \{ u_1 \} $ and $ \kk = (k_1) $. This allows us to write
\begin{align*}
    A_L(f) &= \sum_{u_1 \in \{ 1, \dots, L\}} \sum_{k_1 \leq N_{u_1}} f(\xx_{u_1,k_1}) \sum_{\substack{\widetilde{\fraku} \subseteq \{ u_1 + 1, \dots, L \\ \fraku = \{ u_1 \} \cup \widetilde{\fraku} }} (-1)^{\# \fraku + 1} \cdot \\
    &\phantom{3em} \cdot \sum_{\widetilde{\kk} \leq \NN_{\widetilde{\fraku}}} \left( \prod_{\ell=1}^{\# \fraku - 1} \chi_{u_{\ell},k_{\ell}} ( \xx_{u_{\ell + 1},k_{\ell+1}} ) \right) \chi_{u_{\# \fraku}, k_{\# \kk}}.
\end{align*}
Now we want to rewrite this expression in terms of the points $X_L $ on level $ L $. For a fixed point $ x^* \in X_L $, there is a unique level $ u(\xx^*) $ on which the point occurs first in $ X_{u(\xx^*)} $. By grouping the terms by the points $ \xx \in X_L $ in the representation above, we arrive at
\begin{align*}
    A_L(f) &= \sum_{\xx \in X_L} f(\xx) \sum_{u_1 \in \{u(\xx), \dots, L \}} \sum_{\substack{\widetilde{\fraku} \subseteq \{u_1 + 1, \dots, L\} \\ \fraku = \{u_1\} \cup \widetilde{\fraku}}} (-1)^{\# \fraku + 1} \cdot \\
    &\phantom{3em} \cdot \sum_{\substack{\kk \leq \NN_{\fraku} \\ k_1 = k(\xx,u_1)}} \left( \prod_{\ell=1}^{\# \fraku - 1} \chi_{u_{\ell},k_{\ell}} ( \xx_{u_{\ell + 1},k_{\ell+1}} ) \right) \chi_{u_{\# \fraku}, k_{\# \kk}}.
\end{align*}
Next, we realize that the term 
\begin{align*}
    \sum_{\substack{\kk \leq \NN_{\fraku} \\ k_1 = k(\xx,u_1)}} \left( \prod_{\ell=1}^{\# \fraku - 1} \chi_{u_{\ell},k_{\ell}} ( \xx_{u_{\ell + 1},k_{\ell+1}} ) \right) \chi_{u_{\# \fraku}, k_{\# \kk}}
\end{align*}
is the expression $ \cI_{\fraku} \chi_{L,k(\xx,L)}$, i.e., the combined operator of \eqref{eq:defCombinedOperator} applied to the Lagrange function on level $ L $ centered in the point with index $ k(\xx,L) $ in $ X_L $, we obtain
\begin{align*}
    &\sum_{u_1 \in \{u(\xx), \dots, L \}} \sum_{\substack{\widetilde{\fraku} \subseteq \{u_1 + 1, \dots, L\} \\ \fraku = \{u_1\} \cup \widetilde{\fraku}}} (-1)^{\# \fraku + 1} \cdot \\
    &\phantom{3em} \cdot \sum_{\substack{\kk \leq \NN_{\fraku} \\ k_1 = k(\xx,u_1)}} \left( \prod_{\ell=1}^{\# \fraku - 1} \chi_{u_{\ell},k_{\ell}} ( \xx_{u_{\ell + 1},k_{\ell+1}} ) \right) \chi_{u_{\# \fraku}, k_{\# \kk}} \\
    &= \sum_{u_1 \in \{u(\xx), \dots, L \}} \sum_{\substack{\widetilde{\fraku} \subseteq \{u_1 + 1, \dots, L\} \\ \fraku = \{u_1\} \cup \widetilde{\fraku}}} (-1)^{\# \fraku + 1} \cI_{\fraku} \chi_{L,k(\xx,L)} \\
    &= \sum_{\emptyset \neq \fraku \subseteq \{u(xx), \dots, L} (-1)^{\# \fraku + 1} \cI_{\fraku} \chi_{L,k(\xx,L)} \\
    &=: A_{\{ u(\xx), \dots, L\}} \chi_{L,k(\xx,L)}.
\end{align*}
Hence, we have the representation
\begin{align*}
    A_L(f) = \sum_{\yy \in X_L} f(\yy) A_{\{ u(\yy), \dots, L\}} \chi_{L,k(\yy,L)},
\end{align*}
that was claimed in \cref{thrm:NodalRepresentationMultilevel}.
\end{proof}

\begin{proof}[Proof of \cref{thrm:NodalRepresentationTPML}]
The claim of \cref{thrm:NodalRepresentationTPML} follows from inserting the representation of the direction-wise multilevel operator derived in \cref{thrm:NodalRepresentationMultilevel} into the combination technique representation of the Smolyak operator, together with the observation, that for every component $ \xx^{(j)} $ of elements of the sparse grid $ \cH_{\cI_{\oomega}(\ell,d)} $ occurs first on level $ u_j(\xx^{(j)}) $ and hence, is contained in $ X^{(1)}_{\lambda_1} \times \cdots \times X^{(d)}_{\lambda_d} $ for $ \llambda \in \cI_{\oomega}(\ell,d) $ if and only if $ \lambda_j \geq u_j(\xx^{(j)}) $, $1 \leq j \leq d $. This then yields the claimed nodal representation.
\end{proof}

\newpage
\printbibliography

\end{document}